\documentclass{article}
\usepackage{amsmath,amssymb,pstricks,wasysym,stmaryrd,enumerate,epsfig}

\newcommand{\keywords}[1]{{{\bf Keywords: }#1}}
\newcommand{\tmop}[1]{\operatorname{#1}}
\newcommand{\tmem}[1]{{\em #1\/}}
\newtheorem{definition}{Definition}
\newenvironment{enumerateroman}{\begin{enumerate}[i.]}{\end{enumerate}}
\newcommand{\longrightarrowlim}{\mathop{\longrightarrow}\limits}
\newtheorem{theorem}{Theorem}
\newcommand{\assign}{:=}
\newcommand{\tmstrong}[1]{\textbf{#1}}
\newenvironment{enumeratealpha}{\begin{enumerate}[a{\textup{)}}]}{\end{enumerate}}
\newenvironment{proof}{
  \noindent\textbf{Proof}\ }{\hspace*{\fill}
  \begin{math}\Box\end{math}\medskip}
\newenvironment{proof*}[1]{
  \noindent\textbf{#1\ }}{\hspace*{\fill}
  \begin{math}\Box\end{math}\medskip}
\newcommand{\equallim}{\mathop{=}\limits}
\newtheorem{corollary}{Corollary}
\newtheorem{proposition}{Proposition}
\newcommand{\tmmathbf}[1]{\boldsymbol{#1}}

\begin{document}

\title{Hopf-Galois systems and Kashiwara algebras}\author{Cyril
Grunspan,\\ \\ Institut de Recherche Math\'ematique avanc\'ee,\\
Universit\'e Louis Pasteur et CNRS,\\
7 Rue Ren\'e Descartes,\\ 67084 Strasbourg Cedex, France}\maketitle

\begin{abstract}
  This article is made up with two parts. In the first part, using a recent
  result of Schauenburg, one generalizes to the case when objects are
  faithfully flat over the ground ring, the full equivalence between the
  notions of Hopf-Galois objects and Hopf-Galois systems. In this last
  description, one gives explicitly an inverse for a Hopf-Galois object $T$
  together with its generalized antipode.
  
  In the second part of the article, one shows that the Kashiwara algebras
  introduced by Kashiwara in his study of crystal bases form Hopf-Galois
  systems under the coaction of a quantized enveloping algebra of a Kac-Moody
  algebra. Their classical limits are examples of Sridharan algebras.
  
  {\keywords{Hopf-Galois system, Hopf-Galois extension, Kashiwara algebra,
  Sridharan algebra, torsor.}}
\end{abstract}

\section{Introduction}

In {\cite{Gruns}}, we have introduced the concept of a quantum torsor which
turned out to be an intrinsic reformulation of the old notion of Hopf-Galois
extension {\cite{Schaequi}}, at least in the case when modules are vector
spaces over a field. A quantum torsor is an algebra $T$ equipped with a triple
coproduct $\mu$ satisfying some natural relations of co-associativity and
co-unity. The law $\mu$ ``encodes'' the two structures of bi-algebras of the
Hopf algebras co-acting on $T$. The antipodes are then ``encoded'' themselves
in a torsor endomorphism $\theta$ on $T$ whose existence is a consequence of
torsor axioms satisfied by $\mu$ alone {\cite{qSchau}}. In this last article,
Schauenburg also proved the full equivalence between the two notions of
non-commutative torsors and Hopf-Galois objects to the case when modules are
faithfully flat over the ground ring.

Independently, Bichon introduced in {\cite{Bi}} the concept of Hopf-Galois
systems which is based on the idea that a torsor is a groupoid $( A, B, T, Z
)$ with two units. Originally, this idea can be found in {\cite{Sch}}. With
the help of the Tannaka duality, he proved that if the ground ring is a field,
then the two concepts of Hopf-Galois systems and Hopf-Galois extensions are
equivalent. However, an inverse $Z$ of $T$ is not given explicitly, nor the
``generalized antipode'' from $Z$ to $T$ which is the main feature of his set
of axioms. Note also that his set of axioms is not symmetric. The advantage is
that it is easy to handle with into practice.

In the first part of the article, we show that if $T$ is a $( A, B
)$-bi-Galois object, then the algebra $Z = T^{- 1}$ introduced by Schauenburg
in {\cite{Sch}} can be taken so as to form a Hopf-Galois system $( A, B, T, Z
)$. The result is proved for any faithfully flat Hopf-Galois object $T$. This
generalizes Bichon's theorem. We give explicitly a formula for the
``generalized antipode'' in terms of the torsor endomorphism $\theta$ of $T$.
Moreover, the result one obtains is symmetric in $T$ and $Z$. Note that if $T$
is autonomous i.e., if the torsor map $\theta_{}$ associated to $T$ is
bijective, then one can take $Z = T^{\tmop{op}}$ as proved in {\cite{Gruns}}.

In the second part of this article, we show that the Kashiwara algebras
introduced by Kashiwara in his study of crystal basis {\cite{Kash}} form
Hopf-Galois systems. The Hopf algebras involved in these systems are standard
quantum groups $U_q ( \mathfrak{g}$) associated to any Kac-Moody algebra
$\mathfrak{g}$: Kashiwara algebras are $U_q ( \mathfrak{g} )$-Hopf-Galois
objects. We give the explicit torsor laws $\mu$ and $\theta$. We show that
their classical limits are examples of Sridharan algebras {\cite{Sri}}. 

Kashiwara algebras provide examples of $( A, B )$-bi-Galois objects with two
{\tmem{non-isomorphic}} Hopf-algebras $A$ and $B$ and with a
{\tmem{non-trivial}} torsor map $\theta$, such that at the classical level,
the Hopf algebras associated to the corresponding Hopf-Galois systems are
{\tmem{isomorphic}} and such that the torsor maps are {\tmem{trivial}}. Proofs
of Sections \ref{sridh}, \ref{Kashi} and \ref{class} are very easy and can be
omitted.

\section{Faithfully flat Hopf-Galois systems}

One recalls Bichon's axioms of a Hopf-Galois system {\cite{Bi}}. Fix $k$ a
commutative ring.

\begin{definition}
  \label{defbichon}A Hopf-Galois system consists of four non-zero $k$-module-algebras
  $( A, B, T, Z )$ with the following axioms.
  \begin{enumerateroman}
    \item The algebras $A$ and $B$ are bi-algebras.
    
    \item The algebra $T$ is a $( A, B )$-bicomodule algebra.
    
    \item There are algebra morphisms $\gamma : A \longrightarrow T \otimes Z$
    and $\delta : B \longrightarrow Z \otimes T$ such that the following
    diagrams commute:
    \[ \begin{array}{ccccccccccc}
         T & \longrightarrowlim^{\alpha_T} & A \otimes T &  & A &
         \longrightarrowlim^{\Delta_A} & A \otimes A &  & B &
         \longrightarrowlim^{\Delta_B} & B \otimes B\\
         \beta_T \downarrow &  & \downarrow \gamma \otimes T &  & \gamma
         \downarrow &  & \downarrow A \otimes \gamma &  & \delta \downarrow & 
         & \downarrow \delta \otimes B\\
         T \otimes B & \longrightarrowlim^{T \otimes \delta} & T \otimes Z
         \otimes T &  & T \otimes Z & \longrightarrowlim^{\alpha_T \otimes Z}
         & A \otimes T \otimes Z &  & Z \otimes T & \longrightarrowlim^{Z
         \otimes \beta_T} & Z \otimes T \otimes B
       \end{array} \]
    \item There is a linear map $S_Z : Z \longrightarrow T$ such that the
    following diagrams commute:
    \[ \begin{array}{ccccccccc}
         A & \longrightarrowlim^{\varepsilon_A} k \longrightarrowlim^{u_T} & T
         & ^{} &  &  & B & \longrightarrowlim^{\varepsilon_B} k
         \longrightarrowlim^{u_T} & T\\
         \gamma \downarrow &  & \uparrow m_T &  &  &  & \delta \downarrow &  &
         \uparrow m_T\\
         T \otimes Z & \longrightarrowlim^{T \otimes S_Z} & T \otimes T &  & 
         &  & Z \otimes T & \longrightarrowlim^{S_Z \otimes T} & T \otimes T
       \end{array} \]
  \end{enumerateroman}
\end{definition}

Note that if $( A, B, T, Z )$ is a Hopf-Galois system, then 1) the
bi-algebras $A$ and $B$ are in fact Hopf algebras and 2) the ``generalized
antipode'' $S_Z : Z \longrightarrow T^{\tmop{op}}$ is then an algebra morphism
as proved in {\cite{Bi}}. It is easy to show that if $( A, B, T, Z )$ is a
Hopf-Galois system, then $T$ is in fact a $( A, B )$-bi-Galois extension of
$k$. Conversely, using techniques of Tannaka duality, Bichon proved the
following theorem.

\begin{theorem}
  If $k$ is a field and if $T$ is a $B$-Galois extension of $k$, then there
  exists an Hopf algebra $A$, an algebra $Z$ and a map $S_Z : Z
  \longrightarrow T^{}$ such that $( A, B, T, Z )$ is a Hopf-Galois system.
\end{theorem}

We generalize this theorem to the case when $T$ is a faithfully flat
$B$-Galois extension of $k$. In fact, we prove more. We prove that any
faithfully flat Hopf-Galois objects can be part of a {\tmem{complete}}
Hopf-Galois system in the sense we gave in {\cite{Gruns}} without any
hypothesis on $\theta$.

\begin{theorem}
  \label{bichongene}Let $( T, \mu, \theta )$ be a faithfully flat quantum
  torsor over $k$. Set $A \assign H_l ( T ), B \assign H_r ( T )$ and
  $\alpha_T : T \longrightarrow A \otimes T$ and $\beta_T : T \longrightarrow
  T \otimes B$ the two natural morphisms defining a structure of $( A, B
  )$-bi-Galois extension of $k$ on $T$. Recall that $A$ (resp. $B$) is defined
  as being a certain sub-algebra of $T \otimes T^{\tmop{op}}$ (resp.
  $T^{\tmop{op}} \otimes T$). Then,
  \begin{enumerate}
    \item One has $( B \otimes T )^{\text{co} B} = ( T \otimes A )^{\text{co}
    A}$ as subsets of $T \otimes T \otimes T$.
    
    \item There is a natural structure of $k$-algebra on $Z \assign ( B
    \otimes T )^{\text{co} B} = ( T \otimes A )^{\text{co} A}$ as sub-algebra
    of $T^{\tmop{op}} \otimes T \otimes T^{\tmop{op}}$.
    
    \item There is a natural structure of $( B, A )$-bi-Galois extension of
    $k$ on $Z$ given by morphisms $\beta_Z \assign ( \Delta_B \otimes T ) : Z
    \longrightarrow B \otimes Z$ and $a_Z \assign ( Z \otimes \Delta_A ) : Z
    \longrightarrow Z \otimes A$.
    
    \item If $x \in T$, then $S_T ( x ) \assign ( \theta \otimes T \otimes
    \theta ) \circ \mu^{\tmop{op}} ( x ) \in Z$.
    
    \item The map $S_T : T \longrightarrow Z^{\tmop{op}}$ is an algebra
    morphism.
    
    \item If $x \in Z$, then $( \varepsilon_B \otimes T ) ( x ) = ( T \otimes
    \varepsilon_A ) ( x )$. We denote by $S_Z ( x )$ this common value.
    
    \item The map $S_Z : Z \longrightarrow T^{\tmop{op}}$ is an algebra
    morphism.
    
    \item If $h = x_i \otimes y_i \in B \subset T^{\tmop{op}} \otimes T$, then
    $\delta ( h ) \assign x_i \otimes y^{( 1 )}_i \otimes y_i^{( 2 )} \otimes
    y_i^{( 3 )} \in Z \otimes T \subset T^{\tmop{op}} \otimes T \otimes
    T^{\tmop{op}} \otimes T$.
    
    \item If $h = x_i \otimes y_i \in A \subset T \otimes T^{\tmop{op}}$, then
    $\gamma ( h ) \assign x_i^{( 1 )} \otimes x_i^{( 2 )} \otimes x_i^{( 3 )}
    \otimes y_i \in T \otimes Z \subset T \otimes \text{$T^{\tmop{op}} \otimes
    T \otimes T^{\tmop{op}}$.}$
    
    \item The map $\delta : B \longrightarrow Z \otimes T$ is an algebra
    morphism.
    
    \item The map $\gamma : A \longrightarrow T \otimes Z$ is an algebra
    morphism.
  \end{enumerate}
  Moreover, the quadruple $( A, B, T, Z )$ equipped with the morphisms $( a_T,
  \beta_T, \alpha_Z, \beta_Z, \gamma, \delta, S_T, S_Z$) is a
  {\tmstrong{complete}} Hopf-Galois system in the sense of our definition
  given in {\cite{Gruns}} and the quantum torsor associated to this
  Hopf-Galois system is isomorphic to $( T, \mu, \theta )$. In particular, $(
  A, B, T, Z )$ and $( B, A, Z, T )$ are two Hopf-Galois systems in the sense
  of Definition \ref{defbichon}.
\end{theorem}

Before starting the proof of Theorem \ref{bichongene}, we recall the following
facts taken from {\cite{Gruns}}.
\begin{enumeratealpha}
  \item \label{hopfa}The Hopf algebra $A$ is the set of $x_i \otimes y_i \in T
  \otimes T^{\tmop{op}}$ such that
  \[ x_i^{( 1 )} \otimes x_i^{( 2 )} \otimes \theta ( x_i^{( 3 )} ) \otimes
     y_i = x_i \otimes y_{i^{}}^{( 3 )} \otimes y_i^{( 2 )} \otimes y_i^{( 1
     )} . \]
  \item \label{hopfb}The Hopf algebra $B$ is the set of $x_i \otimes y_i \in
  T^{\tmop{op}} \otimes T^{}$ such that
  \[ x_i \otimes \theta ( y_i^{( 1 )} ) \otimes y_i^{( 2 )} \otimes y_i^{( 3
     )} = x_i^{( 3 )} \otimes x_i^{( 2 )} \otimes x_i^{( 1 )} \otimes y_i . \]
  \item \label{comulta}The comultiplication on $A \subset T \otimes
  T^{\tmop{op}}$ is defined by $\Delta_A ( x_i \otimes y_i ) = x_i^{( 1 )}
  \otimes x_i^{( 2 )} \otimes x_i^{( 3 )} \otimes y_i$.
  
  \item \label{comultb}The comultiplication on $B \subset T^{\tmop{op}}
  \otimes T$ is defined by $\Delta_B ( x_i \otimes y_i ) = x_i \otimes y_i^{(
  1 )} \otimes y_i^{( 2 )} \otimes y_i^{( 3 )}$.
  
  \item \label{antia}The antipode on $A \subset T \otimes T^{\tmop{op}}$ is
  defined by $S_A ( x_i \otimes y_i ) = y_i \otimes \theta ( x_i )$.
  
  \item \label{antib}The antipode on $B \subset T^{\tmop{op}} \otimes T$ is
  defined by $S_B ( x_i \otimes y_i ) = \theta ( y_i ) \otimes x_i$.
  
  \item \label{coua}The counit on $A \subset T \otimes T^{\tmop{op}}$ is
  defined by $\eta_T \circ \varepsilon_A ( x_i \otimes y_i ) = x_i y_i$.
  
  \item \label{coub}The counit on $B \subset T^{\tmop{op}} \otimes T$ is
  defined by $\eta_T \circ \varepsilon_B ( x_i \otimes y_i ) = x_i y_i$.
  
  \item \label{coma}For $x \in T, \alpha_T ( x ) = x^{( 1 )} \otimes x^{( 2 )}
  \otimes x^{( 3 )} \in A \otimes T \subset T \otimes T^{\tmop{op}} \otimes
  T$.
  
  \item \label{comb}For $x \in T, \beta_T ( x ) = x^{( 1 )} \otimes x^{( 2 )}
  \otimes x^{( 3 )} \in T \otimes B \subset T \otimes T^{\tmop{op}} \otimes T$
  
  \item \label{coright}If $H$ is a Hopf algebra and if $M$ and $N$ are two
  $H$-right-comodules, then $( M \otimes N )^{\tmop{co} H}$ is the set of $x_i
  \otimes y_i \in M \otimes N$ such that $x_i \otimes y_{i, N} \otimes y_{i,
  H} = x_{i, M} \otimes y_i^{} \otimes S ( x_{i, H} )$.
  
  \item \label{coleft}Similarly, if $M$ and $N$ are two $H$-left-comodules,
  then $( M \otimes N )^{\tmop{co} A}$ is the set of $x_i \otimes y_i \in M
  \otimes N$ such that $x_{i, H} \otimes x_{i, M} \otimes y_i = S ( y_{i, H} )
  \otimes x_i \otimes y_i$.
\end{enumeratealpha}
\begin{proof}
  In the following, the expression $x \otimes y$ (resp. $x \otimes y \otimes
  z$) will stand for a sum of terms of the form $x_i \otimes y_j \otimes z_k$.
  
  Let $x_{} \otimes y_{} \otimes z_{} \in ( B \otimes T )^{\tmop{co} B}
  \subset T \otimes T \otimes T$. By \ref{comultb}, \ref{antib}, \ref{comb}
  and \ref{coright}, we get :
  \begin{equation}
    x_{} \otimes y_{} \otimes z_{}^{( 1 )} \otimes z_{}^{( 2 )} \otimes
    z_{}^{( 3 )} = x_{} \otimes y_{}^{( 1 )} \otimes z_{} \otimes \theta (
    y_{}^{( 3 )} ) \otimes y_{}^{( 2 )} . \label{defz1}
  \end{equation}
  So,
  \[ x_{} \otimes y_{}^{( 1 )} \otimes y_{}^{( 2 )} \otimes \theta ( y_{}^{(
     3 )} ) \otimes z_{} = x_{} \otimes y_{} \otimes z_{}^{( 3 )} \otimes
     z_{}^{( 2 )} \otimes z_{}^{( 1 )} . \]
  So, by \ref{hopfa} and the fact that $T$ is faithfully flat, we see that
  $x_i \otimes y_j \otimes z_k \in T \otimes A$. Moreover, by \ref{hopfb}, we
  get:
  \begin{equation}
    x_{}^{( 3 )} \otimes x_{}^{( 2 )} \otimes x_{}^{( 1 )} \otimes y_{}
    \otimes z_{} = x_{} \otimes \theta ( y_{}^{( 1 )} ) \otimes y_{}^{( 2 )}
    \otimes y_{}^{( 3 )} \otimes z_{} . \label{defz2}
  \end{equation}
  So,
  \[ x_{}^{( 1 )} \otimes x_{}^{( 2 )} \otimes x_{}^{( 3 )} \otimes y_{}
     \otimes z_{} = y_{}^{( 2 )} \otimes \theta ( y_{}^{( 1 )} ) \otimes x_{}
     \otimes y_{}^{( 3 )} \otimes z_{} . \]
  So, by \ref{comulta}, \ref{antia}, \ref{coma} and \ref{coleft}, we deduce
  that $x_{} \otimes y_{} \otimes z_{} \in ( T \otimes A )^{\tmop{co} A}$.
  Therefore, $( B \otimes T )^{\tmop{co} B} \subset ( T \otimes A )^{\tmop{co}
  A}$ and similarly, $( T \otimes A )^{\tmop{co} A} \subset ( B \otimes T
  )^{\tmop{co} B}$. So, {\tmstrong{1) is proved}}. Assertions {\tmstrong{2)
  and 3) are proved in {\cite{Sch}}}}.
  
  Let $x \in T$. By torsor axioms, we have
  \begin{equation}
    x^{( 1 )} \otimes x^{( 2 )} \otimes \theta ( x^{( 3 )} ) \otimes x^{( 4 )}
    \otimes x^{( 5 )} = x^{( 1 )} \otimes x^{( 2 )^{( 3 )}} \otimes x^{( 2
    )^{( 2 )}} \otimes x^{( 2 )^{( 1 )}} \otimes x^{( 3 )} \label{tortheta}
  \end{equation}
  So,
  \begin{eqnarray*}
    \theta ( x^{( 5 )} ) \otimes \theta ( x^{( 4 )} ) \otimes \theta ( x^{( 3
    )} ) \otimes x^{( 2 )} \otimes \theta ( x^{( 1 )} ) = \theta ( x^{( 3 )} )
    \otimes \theta ( x^{( 2 )} )^{( 1 )} \otimes x^{( 2 )^{( 2 )}} \otimes
    x^{( 2 )^{( 3 )}} \otimes \theta ( x^{( 1 )} ) . &  & 
  \end{eqnarray*}
  So,
  \[ ( \mu^{\tmop{op}} \otimes T \otimes T ) ( \theta ( x^{( 3 )} ) \otimes
     x^{( 2 )} \otimes \theta ( x^{( 1 )} )) = ( T \otimes \theta \otimes T
     \otimes T \otimes T ) \circ ( T \otimes \mu \otimes T ) ( \theta ( x^{( 3
     )} ) \otimes x^{( 2 )} \otimes \theta ( x^{( 1 )} )) \]
  Hence, by \ref{hopfb}, we see that $\theta ( x^{( 3 )} ) \otimes x^{( 2 )}
  \otimes \theta ( x^{( 1 )} ) \in B \otimes T$. On the other hand, by
  (\ref{tortheta}), we get
  \[ \theta ( x^{( 5 )} ) \otimes x^{( 4 )} \otimes \theta ( x^{( 1 )} )
     \otimes \theta ( x^{( 2 )} ) \otimes \theta ( x^{( 3 )} ) = \theta ( x^{(
     3 )} ) \otimes x^{( 2 )^{( 1 )}} \otimes \theta ( x^{( 1 )} ) \otimes
     \theta ( x^{( 2 )^{( 3 )}} ) \otimes x^{( 2 )^{( 2 )}} . \]
  So, by \ref{comultb}, \ref{antib}, \ref{comb} and \ref{coright}, $\theta (
  x^{( 3 )} ) \otimes x^{( 2 )} \otimes \theta ( x^{( 1 )} ) \in ( B \otimes T
  )^{\tmop{co} B}$ and {\tmstrong{4) is proved}}.
  
  For $x, y \in T$, we have $S_T ( y ) \cdot S_T ( x ) = \theta ( x^{( 3 )} )
  \theta ( y^{( 3 )} ) \otimes y^{( 2 )} x^{( 2 )} \otimes \theta ( x^{( 1 )}
  ) \otimes \theta ( y^{( 1 )} )$ since $B \otimes T$ is seen as a subalgebra
  of $T^{\tmop{op}} \otimes T \otimes T^{\tmop{op}}$. Thus, $S_T ( y ) \cdot
  S_T ( x ) =\theta ( x^{( 3 )} y^{( 3 )} ) \otimes y^{( 2 )} x^{( 2 )}
  \otimes \theta ( x^{( 1 )} y^{( 1 )} ) = S_T ( x \cdot y )$. So, $S_T : T
  \longrightarrow Z^{\tmop{op}}$ is an algebra morphism and {\tmstrong{5) is
  proved}}.
  
  Let $u = :x_i \otimes y_j \otimes z_k \in Z$. Then, by \ref{coua} and
  \ref{coub}, $( \varepsilon_B \otimes T ) ( u ) = x_i y_j z_k = ( T \otimes
  \varepsilon_A ) ( u )$. So, {\tmstrong{6) is proved}}.
  
  Suppose that $u = :x_{} \otimes y_{} \otimes z_{} \in Z$ and $v = :x'_{}
  \otimes y'_{} \otimes z'_{} \in Z$. Then, by \ref{coub} and 2),
  \begin{eqnarray*}
    1 \otimes S_Z ( v ) S_Z ( u ) & = & 1 \otimes \varepsilon_B ( x_{}'
    \otimes y' ) z_{}' \cdot \varepsilon_B ( x_{} \otimes y_{} ) z_{} = x'_{}
    y'_{} \otimes z_{}' \cdot \varepsilon_B ( x_{} \otimes y_{} ) z_{}\\
    & = & x_{}' \varepsilon_B ( x_{} \otimes y_{} ) y_{}' \otimes z'_{} z_{}
    = x' _{} x_{} y_{} y_{}' \otimes z'_{} z_{}\\
    & = & 1 \otimes S_Z ( x_{}' x_{} \otimes y_{} y_{}' \otimes z_{}' z_{} )
    = 1 \otimes S_Z ( u v )
  \end{eqnarray*}
  Hence, $S_Z ( u v ) = S_Z ( v u )$ and {\tmstrong{7) is proved}}.
  
  Let $h = x_i \otimes y_i \in B \subset T^{\tmop{op}} \otimes T$. By,
  \ref{comultb} we already have $\delta ( h ) \assign x_{} \otimes y^{( 1
  )}_{} \otimes y_{}^{( 2 )} \otimes y_{}^{( 3 )} \in B \otimes B \subset B
  \otimes T^{} \otimes T$. Moreover, $\delta ( h ) \in ( B \otimes T
  )^{\tmop{co} B} \otimes T$ since by torsor axioms, we have
  \[ x_{} \otimes y_{}^{( 1 )} \otimes y_{}^{( 2 )^{( 1 )}} \otimes y_{}^{( 2
     )^{( 2 )}} \otimes y_{}^{( 2 )^{( 3 )}} \otimes y_{}^{( 3 )} = x_{}
     \otimes y_{}^{( 1 )} \otimes y_{}^{( 4 )} \otimes \theta ( x_{}^{( 3 )} )
     \otimes y_{}^{( 2 )} \otimes y_{}^{( 5 )} . \]
  So, {\tmstrong{8) is proved}}. In the same way, {\tmstrong{9) is true}}.
  
  From the definition of $\delta$ and $\gamma$, it is easy to see that they
  define algebra morphisms. So, {\tmstrong{10) and 11) are true}}.
  
  Let us prove that $( A, B, T, Z )$ is a Hopf-Galois system in the sense of
  Definition \ref{defbichon}. Axioms (i) and (ii) of Definition
  \ref{defbichon} are clearly satisfied. Let $x \in T$. Then,
  \begin{eqnarray*}
    ( \gamma \otimes T ) \circ \alpha_T ( x ) & = & ( \gamma \otimes T ) (
    x^{( 1 )} \otimes x^{( 2 )} \otimes x^{( 3 )} )\\
    & = & x^{( 1 )} \otimes x^{( 2 )} \otimes x^{( 3 )} \otimes x^{( 4 )}
    \otimes x^{( 5 )}\\
    & = & ( T \otimes \delta ) ( x^{( 1 )} \otimes x^{( 2 )} \otimes x^{( 3
    )} )\\
    & = & ( T \otimes \delta ) \circ \beta_T ( x ) .
  \end{eqnarray*}
  This proves the commutativity of the first diagram of (iii) in Definition
  \ref{defbichon}. Let $x \otimes y \in A$. Then,
  \begin{eqnarray*}
    ( A \otimes \gamma ) \circ \Delta_A ( x \otimes y ) & = & ( A \otimes
    \gamma ) ( x^{( 1 )} \otimes x^{( 2 )} \otimes x^{( 3 )} \otimes y )\\
    & = & x^{( 1 )} \otimes x^{( 2 )} \otimes x^{( 3 )} \otimes x^{( 4 )}
    \otimes x^{( 5 )} \otimes y\\
    & = & ( \alpha_T \otimes Z ) ( x^{( 1 )} \otimes x^{( 2 )} \otimes x^{( 3
    )} \otimes y )\\
    & = & ( \alpha_T \otimes Z ) \circ \gamma ( x \otimes y )
  \end{eqnarray*}
  and
  \begin{eqnarray*}
    m_T \circ ( T \otimes S_Z ) \circ \gamma ( x \otimes y ) & = & m_T \circ (
    T \otimes S_Z ) \circ \gamma ( x^{( 1 )} \otimes x^{( 2 )} \otimes x^{( 3
    )} \otimes y )\\
    & = & m_T \circ ( x^{( 1 )} \otimes \varepsilon_B ( x^{( 2 )} \otimes
    x^{( 3 )} ) y )\\
    & = & m_T ( x \otimes y )\\
    & = & \eta_T \circ \varepsilon_A ( x \otimes y )
  \end{eqnarray*}
  Moreover, for $x \otimes y \in B$, we have,
  \begin{eqnarray*}
    ( \delta \otimes B ) \circ \Delta_B ( x \otimes y ) & = & ( \delta \otimes
    B ) ( x \otimes y^{( 1 )} \otimes y^{( 2 )} \otimes y^{( 3 )} )\\
    & = & x \otimes y^{( 1 )} \otimes y^{( 2 )} \otimes y^{( 3 )} \otimes
    y^{( 4 )} \otimes y^{( 5 )}\\
    & = & ( Z \otimes \beta_T ) ( x \otimes y^{( 1 )} \otimes y^{( 2 )}
    \otimes y^{( 3 )} )\\
    & = & ( Z \otimes \beta_T ) \circ \delta ( x \otimes y )
  \end{eqnarray*}
  and
  \begin{eqnarray*}
    m_T \circ ( S_Z \otimes T ) \circ \delta ( x \otimes y ) & = & m_T \circ (
    S_Z \otimes T ) ( x \otimes y^{( 1 )} \otimes y^{( 2 )} \otimes y^{( 3 )}
    )\\
    & = & m_T ( x \varepsilon_A ( y^{( 1 )} \otimes y^{( 2 )} ) \otimes y^{(
    3 )} )\\
    & = & m_T ( x \otimes y )\\
    & = & \eta_T \circ e_B ( x \otimes y ) .
  \end{eqnarray*}
  All these identities above show that {\tmstrong{$( A, B, T, Z )$ is a
  Hopf-Galois system}}.
  
  Let us prove that $( B, A, Z, T )$ equipped with the maps $( \alpha_Z,
  \beta_Z, \ldots, S_T )$ is also a Hopf-Galois system. By 3), Axioms (i) and
  (ii) for this system are satisfied. Let $x \otimes y \otimes z \in Z$. Then,
  \begin{eqnarray*}
    ( Z \otimes \gamma ) \circ \alpha_Z ( x \otimes y \otimes z ) & = & ( Z
    \otimes \gamma ) ( x \otimes y^{( 1 )} \otimes y^{( 2 )} \otimes y^{( 3 )}
    \otimes z )\\
    & = & x \otimes y^{( 1 )} \otimes y^{( 2 )} \otimes y^{( 3 )^{( 1 )}}
    \otimes y^{( 3 )^{( 2 )}} \otimes y^{( 3 )^{( 3 )}} \otimes z\\
    & = & x \otimes y^{( 1 )^{( 1 )}} \otimes y^{( 1 )^{( 2 )}} \otimes y^{(
    1 )^{( 3 )}} \otimes y^{( 2 )} \otimes y^{( 3 )} \otimes z\\
    & = & ( \delta \otimes Z ) \circ \beta_Z ( x \otimes y \otimes z )
  \end{eqnarray*}
  Moreover, for $x \otimes y \in A$, we have
  \begin{eqnarray*}
    ( \gamma \otimes A ) \circ \Delta_A ( x \otimes y ) & = & ( \gamma \otimes
    A ) ( x^{( 1 )} \otimes x^{( 2 )} \otimes x^{( 3 )} \otimes y )\\
    & = & x^{( 1 )^{( 1 )}} \otimes x^{( 1 )^{( 2 )}} \otimes x^{( 1 )^{( 3
    )}} \otimes x^{( 2 )} \otimes x^{( 3 )} \otimes y\\
    & = & x^{( 1 )} \otimes x^{( 2 )} \otimes x^{( 3 )^{( 1 )}} \otimes x^{(
    3 )^{( 2 )}} \otimes x^{( 3 )^{( 3 )}} \otimes y\\
    & = & ( T \otimes \alpha_Z ) ( x^{( 1 )} \otimes x^{( 2 )} \otimes x^{( 3
    )} \otimes y )\\
    & = & ( T \otimes A ) \circ ( T \otimes \alpha_Z ) ( x \otimes y )
  \end{eqnarray*}
  and
  \begin{eqnarray*}
    m_Z \circ ( S_T \otimes Z ) \circ \gamma ( x \otimes y ) & = & m_Z \circ (
    S_T \otimes Z ) ( x^{( 1 )} \otimes x^{( 2 )} \otimes x^{( 3 )} \otimes y
    )\\
    & = & m_Z \circ ( \theta ( x^{( 3 )} ) \otimes x^{( 2 )} \otimes \theta (
    x^{( 1 )} ) \otimes x^{( 4 )} \otimes x^{( 5 )} \otimes y )\\
    & = & x^{( 4 )} \theta ( x^{( 3 )} ) \otimes x^{( 2 )} x^{( 5 )} \otimes
    y \theta ( x^{( 1 )} )\\
    & = & 1 \otimes x^{( 2 )} x^{( 3 )} \otimes y \theta ( x^{( 1 )} )\\
    & = & 1 \otimes 1 \otimes y \theta ( x )\\
    & = & \eta_Z \circ \varepsilon_A ( x \otimes y )
  \end{eqnarray*}
  We also see that for $x \otimes y \in B$, we have
  \begin{eqnarray*}
    ( B \otimes \delta ) \circ \Delta_B ( x \otimes y ) & = & ( B \otimes
    \delta ) ( x \otimes y^{( 1 )} \otimes y^{( 2 )} \otimes y^{( 3 )} )\\
    & = & x \otimes y^{( 1 )} \otimes y^{( 2 )} \otimes y^{( 3 )} \otimes
    y^{( 4 )} \otimes y^{( 5 )}\\
    & = & ( \beta_Z \otimes T ) ( x \otimes y^{( 1 )} \otimes y^{( 2 )}
    \otimes y^{( 3 )} )\\
    & = & ( \beta_Z \otimes T ) \circ \delta ( x \otimes y )
  \end{eqnarray*}
  and
  \begin{eqnarray*}
    m_Z \circ ( Z \otimes S_T ) \circ \delta ( x \otimes y ) & = & m_Z \circ (
    Z \otimes S_T ) ( x \otimes y^{( 1 )} \otimes y^{( 2 )} \otimes y^{( 3 )}
    )\\
    & = & m_Z ( x \otimes y^{( 1 )} \otimes y^{( 2 )} \otimes \theta ( y^{( 5
    )} ) \otimes y^{( 4 )} \otimes \theta ( y^{( 3 )} )\\
    & = & \theta ( y^{( 5 )} ) x \otimes y^{( 1 )} y^{( 4 )} \otimes \theta (
    y^{( 3 )} ) y^{( 2 )}\\
    & = & \theta ( y^{( 3 )} ) x \otimes y^{( 1 )} y^{( 2 )} \otimes 1\\
    & = & \theta ( y ) x \otimes 1 \otimes 1\\
    & = & \eta_Z \circ \varepsilon_B ( x \otimes y )
  \end{eqnarray*}
  We then deduce that {\tmstrong{$( B, A, Z, T )$ is a Hopf-Galois system}}.
  To end the proof, we need to check some compatibility relations between the
  Hopf-Galois systems $( A, B, T, Z )$ and $( B, A, Z, T )$. Let $x \otimes y
  \in A$. Then,
  \begin{eqnarray*}
    ( \beta_T \otimes Z ) \circ \gamma ( x \otimes y ) & = & ( \beta_T \otimes
    Z ) ( x^{( 1 )} \otimes x^{( 2 )} \otimes x^{( 3 )} \otimes y )\\
    & = & x^{( 1 )} \otimes x^{( 2 )} \otimes x^{( 3 )} \otimes x^{( 4 )}
    \otimes x^{( 5 )} \otimes y\\
    & = & ( T \otimes \beta_Z ) ( x^{( 1 )} \otimes x^{( 2 )} \otimes x^{( 3
    )} \otimes y )\\
    & = & ( T \otimes \beta_Z ) \circ \gamma ( x \otimes y )
  \end{eqnarray*}
  and
  \begin{eqnarray*}
    \tau_{( T, Z )} \circ ( S_T \otimes S_Z ) \circ \gamma ( x \otimes y ) & =
    & \tau_{( T, Z )} \circ ( S_T \otimes S_Z ) ( x^{( 1 )} \otimes x^{( 2 )}
    \otimes x^{( 3 )} \otimes y )\\
    & = & \tau_{( T, Z )} ( \theta ( x^{( 3 )} ) \otimes x^{( 2 )} \otimes
    \theta ( x^{( 1 )} ) \otimes \varepsilon_B ( x^{( 4 )} \otimes x^{( 5 )} )
    y )\\
    & = & \tau_{( T, Z )} ( \theta ( x^{( 3 )} ) \otimes x^{( 2 )} \otimes
    \theta ( x^{( 1 )} ) \otimes y )\\
    & = & y \otimes \theta ( x^{( 3 )} ) \otimes x^{( 2 )} \otimes \theta (
    x^{( 1 )} )\\
    & \equallim^{\ref{hopfa}} & y^{( 1 )} \otimes y^{( 2 )} \otimes y^{( 3 )}
    \otimes \theta ( x )\\
    & = & \gamma ( y \otimes \theta ( x )) = \gamma \circ S_A ( x \otimes y )
  \end{eqnarray*}
  Similarly, we prove $( \alpha_Z \otimes T ) \circ \delta = ( Z \otimes
  \alpha_T ) \circ \delta$ and $\tau_{( Z, T )} \circ ( S_Z \otimes S_T )
  \circ \delta = \delta \circ S_B$. Let $x \in T$. Then,
  \begin{eqnarray*}
    \tau_{( A, Z )} \circ ( S_A \otimes S_T ) \circ \alpha_T ( x ) & = &
    \tau_{( A, Z )} \circ ( S_A \otimes S_T ) ( x^{( 1 )} \otimes x^{( 2 )}
    \otimes x^{( 3 )} )\\
    & = & \tau_{( A, Z )} ( x^{( 2 )} \otimes \theta ( x^{( 1 )} ) \otimes
    \theta ( x^{( 5 )} ) \otimes x^{( 4 )} \otimes \theta ( x^{( 3 )} )\\
    & = & \theta ( x^{( 5 )} ) \otimes x^{( 4 )} \otimes \theta ( x^{( 3 )} )
    \otimes x^{( 2 )} \otimes \theta ( x^{( 1 )} )\\
    & = & \theta ( x^{( 3 )} ) \otimes x^{( 2 )^{( 1 )}} \otimes x^{( 2 )^{(
    2 )}} \otimes x^{( 2 )^{( 3 )}} \otimes \theta ( x^{( 1 )} )\\
    & = & \alpha_Z ( \theta ( x^{( 3 )} ) \otimes x^{( 2 )} \otimes \theta (
    x^{( 1 )} ))\\
    & = & \alpha_Z \circ S_T ( x ) .
  \end{eqnarray*}
  In the same way, we prove that $\tau_{( Z, B )} \circ ( S_T \otimes S_B )
  \circ \beta_T = \beta_Z \circ S_T$. Moreover, for $x \otimes y \otimes z \in
  Z$, we have
  \begin{eqnarray*}
    \beta_T \circ S_Z ( x \otimes y \otimes z ) & = & \beta_T ( x \otimes
    \varepsilon_A ( y \otimes z ))\\
    & = & x^{( 1 )} \otimes x^{( 2 )} \otimes x^{( 3 )} \otimes \varepsilon_A
    ( y \otimes z )\\
    & = & y^{( 2 )} \otimes \theta ( y^{( 1 )} ) \otimes x \otimes
    \varepsilon_A ( y^{( 3 )} \otimes z )\\
    & = & y^{( 2 )} \otimes \varepsilon_A ( y^{( 3 )} \otimes z ) \otimes
    \theta ( y^{( 1 )} ) \otimes x\\
    & = & \tau_{( B, T )} \circ ( S_B \otimes S_Z ) ( x \otimes y^{( 1 )}
    \otimes y^{( 2 )} \otimes y^{( 3 )} \otimes z )\\
    & = & \tau_{( B, T )} \circ ( S_B \otimes S_Z ) \circ \beta_Z ( x \otimes
    y \otimes z )
  \end{eqnarray*}
  In the same way, we prove that $\alpha_T \circ S_Z = \tau_{( T, A )} \circ (
  S_Z \otimes S_A ) \circ \alpha_Z$. Therefore, all the axioms of
  compatibility between the two Hopf-Galois systems $( A, B, T, Z )$ and $( B,
  A, Z, T )$ are satisfied (see {\cite{Gruns}}). Thus, {\tmstrong{$( A, B, T,
  Z )$ is a complete Hopf-Galois system}}.

\end{proof}

Given the full equivalence between Hopf-Galois extensions and quantum torsors
as proved in {\cite{qSchau}}, we can state the following generalization of
Bichon's theorem.

\begin{corollary}
  If $T$ is a faithfully flat $B$-Galois extension of $k$, then there exists
  an Hopf algebra $A$, an algebra $Z$ and algebra morphisms $S_T : T
  \longrightarrow Z^{\tmop{op}}$ and $S_Z : Z \longrightarrow T^{\tmop{op}}$
  such that $( A, B, T, Z )$ is a complete Hopf-Galois system in the sense we
  gave in {\cite{Gruns}}. In particular, $( A, B, T, Z )$ and $( B, A, Z, T )$
  are two Hopf-Galois systems in the sense of Definition \ref{defbichon}.
\end{corollary}

\section{\label{sridh}Sridharan algebras and $\tmop{Bigal} ( U ( \mathfrak{g}
))$}

In this short section, we introduce the Sridharan algebras {\cite{Sri}} and we
prove that they are $U ( \mathfrak{g} )$-Galois objects.

\begin{definition}
  Let $( \mathfrak{g}, [, ] )$ be a $\mathbb{Q}$-Lie algebra and let $c$ be a
  $2$-cocycle of $\mathfrak{g}$ with values in $\mathbb{Q}$ i.e., $c$
  satisfies the following equation:
  \[ c ( [ x, y ], z ) + c ( [ y, z ], x ) + c ( [ z, x ], y ) = 0 \]
  for all $x, y, z \in \mathfrak{g}$. The Sridharan algebra $U_c (
  \mathfrak{g} )$ associated to $c$ is $T ( \mathfrak{g} ) / I_c (
  \mathfrak{g} )$ where $T ( \mathfrak{g} )$ is the tensorial algebra of
  $\mathfrak{g}$ and $I_c ( \mathfrak{g} )$ is the ideal generated by elements
  $x \cdot y - y \cdot x - [ x, y ] - c ( x, y ) \cdot 1$ with $x, y \in
  \mathfrak{g}$.
\end{definition}

It can be shown that two cocycles define two isomorphic Sridharan algebras if
and only if the two cocycles are cohomologous. The point is that Sridharan are
examples of $U ( \mathfrak{g} )$-torsors.

\begin{theorem}
  Let $\mathfrak{g}$ be a $\mathbb{Q}$-Lie algebra and $c$ a $2$-cocycle of
  $\mathfrak{g}$ with values in $\mathbb{Q}$. Then, $( U_c ( \mathfrak{g} ),
  U_{- c} ( \mathfrak{g} ))$ is a $( U ( \mathfrak{g} ), U ( \mathfrak{g}
  ))$-Hopf-Galois system over $\mathbb{Q}$. The maps defining this Hopf-Galois
  system are all of the form $x \mapsto 1 \otimes x + x \otimes 1$ for $x \in
  \mathfrak{g}$ and $S ( x ) = - x, x \in \mathfrak{g}$ for the generalized
  antipodes.
\end{theorem}

By {\cite{Gruns}}, we then deduce the following corollaries.

\begin{corollary}
  Let $\mathfrak{g}$ be a $\mathbb{Q}$-Lie algebra and $c$ a $2$-cocycle of
  $\mathfrak{g}$ with values in $\mathbb{Q}$. The Sridharan algebra $U_c (
  \mathfrak{g} )$ is a quantum torsor and $H_l ( U_c ( \mathfrak{g} )) \cong
  H_r ( U_c ( \mathfrak{g} )) \cong U ( \mathfrak{g} )$. The torsor laws are
  defined by $\theta_{U_c ( \mathfrak{g} )} = \tmop{Id}_{U_c ( \mathfrak{g}
  )}$ and $\mu ( x ) = x \otimes 1 \otimes 1 - 1 \otimes x \otimes 1 + 1
  \otimes 1 \otimes x$ for $x \in \mathfrak{g}$.
\end{corollary}

\begin{corollary}
  Sridharan algebras are $( U ( \mathfrak{g} ), U ( \mathfrak{g}
  ))$-Hopf-bi-Galois extensions of $\mathbb{Q}$.
\end{corollary}

\begin{proposition}
  Let $\mathfrak{g}$ be a $\mathbb{Q}$-Lie algebra and $c, c'$ two
  $2$-cocycles of $\mathfrak{g}$ with values in $\mathbb{Q}$. Then, $U_c (
  \mathfrak{g} ) \Box_{U ( \mathfrak{g} )} U_{c'} ( \mathfrak{g} ) \cong U_{c
  \cdot c'} ( \mathfrak{g} )$, where $A \Box_H B$ denotes the $H$-coproduct of
  a $H$-right-comodule-algebra $A$ and a $H$-left-comodule-algebra $B$. The
  isomorphism is an isomorphism of $( U ( \mathfrak{g} ), U ( \mathfrak{g}
  ))$-Hopf-bi-Galois extensions of $\mathbb{Q}$.
\end{proposition}

\begin{corollary}
  The map $c \longrightarrow U_c ( \mathfrak{g} )$ defines a monomorphism of
  groups of $H^2 ( \mathfrak{g} )$ into $\tmop{Bigal} ( U ( \mathfrak{g} ))$.
\end{corollary}

We conjecture that $\tmop{Bigal} ( U ( \mathfrak{g} )) = H^2 ( \mathfrak{g}
)$. This result is may be known of specialists but we are unable to give a
reference. 

\section{\label{Kashi}Kashiwara algebras}

In Section \ref{Kashqt}, we recall the definition of Kashiwara algebras and we
prove that they can be equipped with a quantum torsor structure. Then, in
Section \ref{Kashhgs}, we will identify the two Hopf algebras associated with
these quantum torsors. First, we fix some notations.

\subsection{Notations.}{\tmstrong{}} We fix $\mathfrak{g}$ a symmetrizable
Kac-Moody algebra over $\mathbb{Q}$ with a Cartan sub-algebra $\mathfrak{t}$,
$\{ \alpha_i \in \mathfrak{t}^{\ast} \}_{i \in I}$ the set of simple roots and
$\{ h_i \in \mathfrak{t} \}_{i \in I}$ the set of coroots, where $I$ is a
finite index set. We define an inner product on $\mathfrak{t}^{\ast}$ such
that $( \alpha_i, \alpha_i ) \in \mathbb{N}$ and $\langle h_i, \lambda \rangle
= 2 ( \alpha_i, \lambda ) / ( \alpha_i, \alpha_i )$ for $\lambda \in
t^{\ast}$. Set $Q = \oplus_i \mathbb{Z} \alpha_i$, $P \equallim \oplus_i
\mathbb{Z} h_i$ and $P^{\ast} = \{ h \in \mathfrak{t} / \langle h, P \rangle
\subset \mathbb{Z} \}$. We set $q = \exp ( \hbar ), q_i = q^{\langle \alpha_i,
\alpha_i \rangle / 2}, t_i = q^{h_i}, [ n ]_i = ( q_i^n - q_i^{- n} ) / ( q_i
- q_i^{- 1} ), [ n ]_i ! = \prod_{k = 1}^n [ k ]_i$. The ground ring $R$ is
$\mathbb{Q} [[ \hbar ]]$. We denote by $F$ the field of Laurent series
$\mathbb{Q} (( \hbar ))$. If $M$ is a $R$-module, we will set $M_F \assign M
\otimes_R F$.

\subsection{\label{Kashqt}Kashiwara algebras and quantum torsors}

Usually, Kashiwara algebras are defined over $F$. However, it can be useful
to see them as topologically free $R$-modules, especially because we will
examine their classical limits.

\begin{definition}
  The Kashiwara algebra $B_q ( \mathfrak{g} )$ {\cite{Ka}} is the associative
  $R$-algebra generated by generators $e'_i, f_i, i \in I, q^h, h \in P$ and
  relations :
  \begin{eqnarray}
    q^h e'_i q^{- h} & = & q^{\langle h, \alpha_i \rangle} e'_i \\
    q^h f_i q^{- h} & = & q^{- \langle h, \alpha_i \rangle} f_i \\
    e'_i f_j & = & q_i^{\langle h_i, \alpha_j \rangle} f_j e'_i + \delta_{i,
    j} . \\
    \sum_{k = 0}^{1 - \langle h_i, \alpha_j \rangle} ( - 1 )^k X_i^{( k )} X_j
    X_i^{( 1 - \langle h_i, \alpha_j \rangle )} & = & 0 \nonumber
  \end{eqnarray}
  with $X = e', f$ and $X_i^{( n )} = X_i^n / [ n ]_i !$. We denote in brief
  $B$ for $B_q ( \mathfrak{g} )$.
\end{definition}

\begin{proposition}
  The algebra $B$ is a topologically free $R$-module. A basis for $B$ over $R$
  is given by the family $(e'_{i_1})^{\alpha} q^h f_{i_2}^{\beta}$
  with $\alpha, \beta \in \mathbb{N}$, $i_1, i_2 \in I$ and $h \in P$.
\end{proposition}

The following theorem proves the existence of a quantum torsor structure on
$B$.

\begin{theorem}
  \label{bqt}The morphisms $\mu_B$ and $\theta_B$ below are well defined and
  $( B, \mu_B, \theta_B )$ is an autonomous quantum torsor {\cite{Gruns}}.
  \[ \begin{array}{cccc}
       \mu_{B_{}} : & B_{} & \longrightarrow & B_{} \otimes B_{}^{\tmop{op}}
       \otimes B\\
       & e'_i & \longmapsto & 1 \otimes 1 \otimes e'_i - 1 \otimes e'_i t_i
       \otimes t_i^{- 1} + e'_i \otimes t_i \otimes t_i^{- 1}\\
       & f_i & \longmapsto & 1 \otimes 1 \otimes f_i - 1 \otimes f_i t_i
       \otimes t_i^{- 1} + f_i \otimes t_i \otimes t_i^{- 1}\\
       & q^h & \longmapsto & q^h \otimes q^{- h} \otimes q^h\\
       &  &  & \\
       \theta_{B_{}} : & B_{} & \longrightarrow & B_{}\\
       & e'_i & \longmapsto & t_i^{- 1} e'_i t_i\\
       & f_i & \longmapsto & t_i^{- 1} f_i t_i\\
       & q^h & \longmapsto & q^h
     \end{array} \]
\end{theorem}

We want now to identify the two Hopf algebras $H_l ( B )$ and $H_r ( B )$
associated with the quantum torsor $( B, \mu_B, \theta_B )$. To this end, we
need to define some quantum groups.

\subsection{\label{defqg}The quantum groups $U_q ( \mathfrak{g} ), U'_q (
\mathfrak{g} )$ and $\hat{U}_q ( \mathfrak{g} )$}

In this section, we will define three quantum groups which will turn out to
coact on $B$. The first one is the Drinfeld-Jimbo standard quantum group.

\begin{definition}
  The quantum group $U \assign U_q ( \mathfrak{g} )$ is the associative
  $R$-algebra over $R$ generated by $\tmmathbf{e}_i, \tmmathbf{f}_i,
  \tmmathbf{q}^h$ ($i \in I, h \in P^{\ast} )$ submitted to the following
  relations:
  \begin{eqnarray}
    \tmmathbf{q^h}  \tmmathbf{e_i}  \tmmathbf{q^{- h}} & = &
    \tmmathbf{q^{\langle h, \alpha_i \rangle}}  \tmmathbf{e_i} \\
    \tmmathbf{q^h}  \tmmathbf{f_i}  \tmmathbf{q^{- h}} & = & \tmmathbf{q^{-
    \langle h, \alpha_i \rangle}}  \tmmathbf{f_i} \\
    \tmmathbf{[ e_i, f_j ]} & = & \delta_{i, j}  \tmmathbf{( t_i - t_i^{- 1}
    )} / ( q_i - q_i^{- 1} ) \linebreak \\
    \sum_{k = 0}^{1 - \langle h_i, \alpha_j \rangle} ( - 1 )^k 
    \tmmathbf{X_i^{( k )}} \tmmathbf{X_j X_i^{( 1 - \langle h_i, \alpha_j
    \rangle )}} & = & 0 
  \end{eqnarray}
  with $\tmmathbf{X = e, f}$ and $\tmmathbf{X_i^{( n )}} = \tmmathbf{X_i^n} /
  [ n ]_i !$.
  
  The co-multiplication $\Delta_{}$ on $U$ is the algebra homomorphism defined
  by:

  \begin{eqnarray}
    \Delta_{} ( \tmmathbf{q}^h ) & = & \tmmathbf{q}^h \otimes \tmmathbf{q}^h
    \\
    \Delta_{} ( \tmmathbf{e}_i ) & = & \tmmathbf{e}_i \otimes 1 +
    \tmmathbf{t}_i \otimes \tmmathbf{e}_i \\
    \Delta_{} ( \tmmathbf{f}_i ) & = & \tmmathbf{f}_i \otimes
    \tmmathbf{t}_i^{- 1} + 1 \otimes \tmmathbf{f}_i . 
  \end{eqnarray}
  The antipode $S : U_{} \longrightarrow U_{} $is the anti-isomorphism defined
  by the formulas:
  \begin{eqnarray}
    S ( \tmmathbf{e}_i ) & = & - \tmmathbf{t}_i^{- 1} \tmmathbf{e}_i \\
    S ( \tmmathbf{f}_i ) & = & - \tmmathbf{f}_i \\
    S ( \tmmathbf{q}^h ) & = & \tmmathbf{q}^{- h} . 
  \end{eqnarray}
  The counit $\varepsilon : U_{} \longrightarrow R$ is the algebra
  homomorphism defined by:
  \begin{eqnarray}
    \varepsilon ( \tmmathbf{e}_i ) & = & 0 \\
    \varepsilon ( \tmmathbf{f}_i ) & = & 0 \\
    \varepsilon ( \tmmathbf{q}^h ) & = & 1 . 
  \end{eqnarray}
\end{definition}

We define now the two quantum groups $U'_q ( \mathfrak{g} )$ and $\hat{U}_q (
\mathfrak{g} )$. They correspond to two different quantizations of a common
enveloping algebra $U ( \mathfrak{a} )$ where $\mathfrak{a}$ is the Lie
algebra defined below.

\begin{definition}
  Let $\mathfrak{b}_+$ and $\mathfrak{b}_-$ be the two Borel parts of
  $\mathfrak{g}$. We denote by $\mathfrak{a}$ the sub-Lie algebra of pairs $(
  x_+, x_- ) \in \mathfrak{b}_+ \times \mathfrak{b}_-$ such that there is $h
  \in \mathfrak{t}$ satisfying $x_{\pm} - h \in \mathfrak{n}_{\pm}$.
\end{definition}

\begin{definition}
  The quantum group $\hat{U}_q ( \mathfrak{g} )$ (denoted in brief by
  $\hat{U}$) is the associative algebra generated by generators
  $\tmmathbf{\hat{e}_i}, \tmmathbf{f_i,} i \in I$ and $\tmmathbf{q^h}, h \in
  P^{\ast}$ satisfying
  \begin{eqnarray}
    \tmmathbf{q^h \tmmathbf{\widehat{e}_i}}  \tmmathbf{} \tmmathbf{q^{- h}} &
    = & \tmmathbf{q^{\langle h, \alpha_i \rangle} \tmmathbf{\hat{e}_i}}  \\
    \tmmathbf{q^h}  \tmmathbf{f_i}  \tmmathbf{q^{- h}} & = & \tmmathbf{q^{-
    \langle h, \alpha_i \rangle}}  \tmmathbf{f_i} \\
    \tmmathbf{[ \hat{e}_i, f_j ]} & = & \delta_{i, j}  \tmmathbf{( t_i -
    t_i^{- 1} )} \linebreak \\
    \sum_{k = 0}^{1 - \langle h_i, \alpha_j \rangle} ( - 1 )^k 
    \tmmathbf{X_i^{( k )}} \tmmathbf{X_j X_i^{( 1 - \langle h_i, \alpha_j
    \rangle )}} & = & 0 
  \end{eqnarray}
  with $\tmmathbf{X = e, f}$ and $\tmmathbf{X_i^{( n )}} = \tmmathbf{X_i^n} /
  [ n ]_i !$.
  
  The co-multiplication $\hat{\Delta}$ on $\hat{U}$ is the algebra
  homomorphism defined by :
  \begin{eqnarray}
    \text{$\hat{\Delta}$}_{} ( \tmmathbf{q}^h ) & = & \tmmathbf{q}^h \otimes
    \tmmathbf{q}^h \\
    \text{$\hat{\Delta}$}_{} ( \tmmathbf{\hat{e}_i} ) & = &
    \tmmathbf{\hat{e}_i} \otimes 1 + \tmmathbf{t}_i \otimes
    \tmmathbf{\hat{e}_i} \\
    \text{$\hat{\Delta}$}_{} ( \tmmathbf{f}_i ) & = & \tmmathbf{f}_i \otimes
    \tmmathbf{t}_i^{- 1} + 1 \otimes \tmmathbf{f}_i . 
  \end{eqnarray}
  The antipode $\hat{S} : \hat{U} \longrightarrow \hat{U}$ is the
  anti-isomorphism defined by the formulas:
  \begin{eqnarray}
    \hat{S} ( \tmmathbf{\hat{e}_i} \tmmathbf{} ) & = & - \tmmathbf{t}_i^{- 1}
    \tmmathbf{\hat{e}_i} \\
    \hat{S} ( \tmmathbf{f}_i ) & = & - \tmmathbf{f}_i \\
    \hat{S} ( \tmmathbf{q}^h ) & = & \tmmathbf{q}^{- h} . 
  \end{eqnarray}
  The counit $ \hat{\varepsilon} : \hat{U} \longrightarrow R$ is the algebra
  homomorphism defined by:
  \begin{eqnarray}
    \hat{\varepsilon} ( \tmmathbf{\hat{e}_i} ) & = & 0 \\
    \hat{\varepsilon} ( \tmmathbf{f}_i ) & = & 0 \\
    \hat{\varepsilon} ( \tmmathbf{q}^h ) & = & 1 . 
  \end{eqnarray}
\end{definition}

It can be shown that the above relations defines effectively a Hopf algebra.
{\tmem{A priori}}, the common notation $\tmmathbf{f}_i$ for the generators of
$U$ and $\hat{U}$ could be ambiguous for further computations. It is not the
case as expressed by the following proposition which links $\hat{U}$ with $U$.

\begin{proposition}
  The morphism from $\hat{U}$ to $U$ which sends $\tmmathbf{}
  \tmmathbf{\hat{e}_i}, \tmmathbf{f_i}, \tmmathbf{q^h}$ to $( q_i - q_i^{- 1}
  ) \tmmathbf{e_i}, \tmmathbf{f_i}, \tmmathbf{q^h}$ is a monomorphism of
  Hopf-algebras.
\end{proposition}

From this, we deduce easily the following corollary.

\begin{corollary}
  \label{form}The quantum group $\hat{U}$ is a $F$-form of $U : \hat{U}
  \otimes F = U \otimes F$
\end{corollary}

Next, we turn to $U'_q ( \mathfrak{g} )$.

\begin{definition}
  The quantum group $U'_q ( \mathfrak{g} )$ (denoted in brief by $U'$) is the
  associative $R$-algebra over $R$ generated by generators $\tmmathbf{e}'_i,
  \tmmathbf{f}'_i, i \in I$ and $\tmmathbf{q}^h, h \in P^{\ast}$ submitted to
  the following relations:
  \begin{eqnarray*}
    \tmmathbf{q^h e'_i q^{- h}} & = & \tmmathbf{q^{\langle h, \alpha_i
    \rangle} e'_i}\\
    \tmmathbf{q^h f'_i q^{- h}} & = & \tmmathbf{q^{- \langle h, \alpha_i
    \rangle} f'_i}\\
    \tmmathbf{e'_i f'_j} & = & \tmmathbf{q_i^{- \langle h_i, \alpha_j \rangle}
    f'_j e'_i}\\
    \sum_{k = 0}^{1 - \langle h_i, \alpha_j \rangle} ( - 1 )^k 
    \tmmathbf{X_i^{( k )}} \tmmathbf{X_j X_i^{( 1 - \langle h_i, \alpha_j
    \rangle )}} & = & 0
  \end{eqnarray*}
  with $\tmmathbf{X = e', f'}$ and $\tmmathbf{X_i^{( n )}} = \tmmathbf{X_i^n}
  / [ n ]_i !$. The co-multiplication $\Delta_{}'$ on $U'$ is the algebra
  homomorphism defined by :
  \begin{eqnarray}
    \Delta' ( \tmmathbf{q}^h ) & = & \tmmathbf{q}^h \otimes \tmmathbf{q}^h \\
    \Delta_{}' ( \tmmathbf{e}'_i ) & = & \tmmathbf{e}'_i \otimes 1 +
    \tmmathbf{t}_i \otimes \tmmathbf{e}'_i \\
    \Delta_{}' ( \tmmathbf{f}'_i ) & = & \tmmathbf{f}'_i \otimes 1 +
    \tmmathbf{t}_i \otimes \tmmathbf{f}'_i . 
  \end{eqnarray}
  The antipode $S' : U_{}' \longrightarrow U_{}' $is the anti-isomorphism
  defined by the formulas:
  \begin{eqnarray}
    S' ( \tmmathbf{e}'_i ) & = & - \tmmathbf{t}_i^{- 1} \tmmathbf{e}'_i \\
    S' ( \tmmathbf{f}'_i ) & = & \tmmathbf{t}_i^{- 1} \tmmathbf{f}'_i \\
    S' ( \tmmathbf{q}^h ) & = & \tmmathbf{q}^{- h} . 
  \end{eqnarray}
  The counit $\varepsilon' : U_{}' \longrightarrow R$ is the algebra
  homomorphism defined by:
  \begin{eqnarray}
    \varepsilon' ( \tmmathbf{e}'_i ) & = & 0 \\
    \varepsilon' ( \tmmathbf{f}'_i ) & = & 0 \\
    \varepsilon' ( \tmmathbf{q}^h ) & = & 1 . 
  \end{eqnarray}
\end{definition}

Here also, it can be shown that the above relations defines effectively a Hopf
algebra. We are now ready to identify $H_l ( B )$ and $H_r ( B )$.

\subsection{\label{Kashhgs}Kashiwara algebras and Hopf-Galois systems}

\subsubsection{Hopf Galois systems and Kashiwara algebras over $R$}

We are going to embed $( B, B^{\tmop{op}} )$ into a structure of autonomous
Hopf-Galois systems.

\begin{proposition}
  The morphisms $\alpha_B, \beta_B, \alpha_{B^{\tmop{op}}},
  \beta_{B^{\tmop{op}}}, \gamma, \delta, S_B, S_{B^{\tmop{op}}}$ below are
  well defined and $( U', \hat{U}, B, B^{\tmop{op}} )$ is a (complete)
  autonomous Hopf-Galois system.
  \[ \begin{array}{cccc}
       \alpha_B : & B & \longrightarrow & U_{}' \otimes B_{}\\
       & e'_i & \longmapsto & 1 \otimes e'_i + \tmmathbf{t_i^{- 1} e'_i}
       \otimes t_i^{- 1}\\
       & f_i & \longmapsto & 1 \otimes f_i + \tmmathbf{t_i^{- 1} f'_i}
       \otimes t_i^{- 1}\\
       & q^h & \longmapsto & \tmmathbf{q^h} \otimes q^h \tmmathbf{}\\
       &  &  & \\
       \beta_B : & B_{} & \longrightarrow & B \otimes \hat{U}\\
       & e'_i & \longmapsto & 1 \otimes \tmmathbf{t_i^{- 1} \hat{e}_i} + e'_i
       \otimes \tmmathbf{t_i^{- 1}}\\
       & f_i & \longmapsto & f_i \otimes \tmmathbf{t_i^{- 1}} + 1 \otimes
       \tmmathbf{f_i}\\
       & q^h & \longmapsto & q^h \otimes \tmmathbf{q^h}\\
       &  &  & \\
       \alpha_{B^{\tmop{op}}} : & B^{\tmop{op}} & \longrightarrow & \hat{U}
       \otimes B^{\tmop{op}}\\
       & e'_i & \longmapsto & \tmmathbf{t_i} \otimes e'_i -
       \tmmathbf{\hat{e}_i} \otimes 1\\
       & f_i & \longmapsto & - \tmmathbf{t_i f_i} \otimes 1 + \tmmathbf{t_i}
       \otimes f_i\\
       & q^h & \longmapsto & \tmmathbf{q^{- h}} \otimes q^h\\
       &  &  & \\
       \beta_{B^{\tmop{op}}} : & B^{\tmop{op}} & \longrightarrow &
       B^{\tmop{op}}_{} \otimes U'\\
       & e'_i & \longmapsto & - t^{- 1}_i \otimes \tmmathbf{e'_i} + e'_i
       \otimes 1\\
       & f_i & \longmapsto & - t^{- 1}_i \otimes \tmmathbf{f'_i} + f_i
       \otimes 1\\
       & q^h & \longmapsto & q^h \otimes \tmmathbf{q^{- h}}\\
       &  &  & \\
       \gamma : & U_{}' & \longrightarrow & B_{} \otimes B^{\tmop{op}}\\
       & \tmmathbf{e'_i} & \longmapsto & t_i e'_i \otimes 1 - t_i \otimes
       e'_i\\
       & \tmmathbf{f'_i} & \longmapsto & t_i f_i \otimes 1 - t_i \otimes
       f_i\\
       & \tmmathbf{q^h} & \longmapsto & q^h \otimes q^{- h}\\
       &  &  & \\
       \delta : & \hat{U} & \longrightarrow & B^{\tmop{op}} \otimes B\\
       & \tmmathbf{\hat{e}_i} & \longmapsto & t_i^{- 1} \otimes t_i e'_i -
       e'_i \otimes 1\\
       & \tmmathbf{f_i} & \longmapsto & 1 \otimes f_i - f_i t_i \otimes
       t_i^{- 1}\\
       & \tmmathbf{q^h} & \longmapsto & q^{- h} \otimes q^h\\
       &  &  & \\
       S_{B^{\tmop{op}}} \assign \tmop{Id}_B : & B^{\tmop{op}} &
       \longrightarrow & ( B )^{\tmop{op}}\\
       &  &  & \\
       S_{B} \assign \theta_B : & B & \longrightarrow & ( B^{\tmop{op}}
       )^{\tmop{op}} = B
     \end{array} \]
\end{proposition}

With the help of {\cite{Gruns}}, we can now identify $H_l ( B )$ and $H_r ( B
)$.

\begin{corollary}
  The morphisms $\gamma$ and $\delta$ define Hopf algebras isomorphisms
  between $U'$ (resp. $\hat{U}$) and $H_l ( B )$ (resp. $H_r ( B )$).
\end{corollary}

We are now going to extend the scalars to $F$.

\subsubsection{Hopf Galois systems and Kashiwara algebras over $F$}

The quadruple $( B, B^{\tmop{op}}, U', \hat{U} )$ is a complete autonomous
Hopf-Galois system. So by {\cite{Gruns}}, we see that $( B_F, B_F^{\tmop{op}},
U'_F, \hat{U}_F )$ is also a complete autonomous Hopf-Galois system. It
follows that $H_r ( B_F ) \cong \hat{U}_F$ and $H_l ( B_F ) \cong U'_F$.
Therefore, with the help of Corollary \ref{form}, we can state the following
theorem.

\begin{theorem}
  Kashiwara algebras equipped with their natural structures of $U_{\hbar} (
  \mathfrak{g} )$-comodules are $U_{\hbar} ( \mathfrak{g} )$-Galois extensions
  of $\mathbb{Q} (( \hbar ))$.
\end{theorem}

\section{\label{class}Classical limit of $B$}

If we consider the limit when $\hbar \mapsto 0$, then $\hat{e}_i$ and $f_j$
commute (resp. $e'_i$ and $f'_j$) and the only relations between generators
are Serre relations. The classical limit of $\hat{U}$ and the classical limit
of $U'$ coincide coincide both with the enveloping algebra $U ( \mathfrak{a}
)$ ($\mathfrak{a}$ has been introduced in the first section). The classical
limit $B_{\tmop{cl}}$ of $B$ is the algebra generated by $e'_i, f_j$ ($i, j
\in I$) and $h \in \mathfrak{t}$ submitted to relations:
\begin{eqnarray}
  [ h, e'_i ] & = & \langle \alpha_i, h \rangle e'_i \label{hepcl} \\
  {[} h, f_j ] & = & - \langle \alpha_i, h \rangle f_j \label{hfcl} \\
  e'_i f_j & = & f_j e'_i + \delta_{i, j} \label{epfcl} 
\end{eqnarray}
together with Serre relations. To the limit when $\hbar \mapsto 0$, we see
that the classical limit of Kashiwara algebras $B_{\tmop{cl}}$ are $( U (
\mathfrak{a} ), U ( \mathfrak{a} ))$-quantum torsors over $\mathbb{Q} .$ In
fact, it is easy to see that classical limit of Kashiwara algebras are
examples of Sridharan algebras.

\begin{proposition}
  The classical limit of Kashiwara algebras $B_{\tmop{cl}}$ are Sridharan
  algebras for the Lie algebra $\mathfrak{a}$ defined above. The cocycle $c :
  \mathfrak{a} \wedge \mathfrak{a} \longrightarrow \mathfrak{a}$ is defined by
  $c ( e_i, f_j ) = \delta_{i, j}, c ( e_i, e_j ) = c ( f_i, f_j ) = c ( h,
  e_i ) = c ( h, f_i ) = 0$ and cocycle relations.
\end{proposition}

\section{Open problems}

We finish this article by asking two questions.

\subsubsection{Question 1. }If $T$ is a $B$-Galois extension of $k$ and 
if the antipode $S_B$ of $B$ is bijective, does this imply $\theta_T$ bijective ?

\subsubsection{Question 2 (Completion problem).} If $( A, B, T, Z )$ is a
Hopf-Galois system in the sense of Bichon, is it possible to find a $( B, A
)$-bicomodule structure on $Z$ and a map $S_T : T \longrightarrow Z$ such that
$( A, B, T, Z )$ is a {\tmem{complete}} Hopf-Galois system in the sense we
gave in {\cite{Gruns}} ?

\subsection{Acknowledgment} I am very grateful to B. Enriquez for his help and
support. I am also grateful to J. Bichon and P. Schauenburg for remarks. I
would also like to thank W. Soergel for the hospitality of the Freiburg
University.

\bigskip

{}

\end{document}